\documentclass{amsproc}
\usepackage{epsfig}
\newtheorem{theorem}{Theorem}[section]

\theoremstyle{definition}
\newtheorem{definition}[theorem]{Definition}

\newtheorem{example}[theorem]{Example}

\theoremstyle{remark}
\newtheorem{remark}[theorem]{Remark}

\numberwithin{equation}{section}

\begin{document}

\title{Secondary Cohomology and $k$-invariants}
\author{Mihai D. Staic}
\address{Department of Mathematics, Indiana University, Rawles Hall, Bloomington, IN 47405, USA }
\thanks{\textit{Permanent address:} Institute of Mathematics of the Romanian Academy, PO.BOX 1-764, RO-70700 Bucharest, Romania. Research partially supported by the
CNCSIS project ''Hopf algebras, cyclic homology and monoidal categories''
contract nr. 560/2009.}
\email{mstaic@indiana.edu}



\date{January 1, 1994 and, in revised form, June 22, 1994.}



\begin{abstract} For a triple $(G,A,\kappa)$ (where $G$ is a group, $A$ is a $G$-module and $\kappa:G^3\to A$ is a 3-cocycle) and a $G$-module $B$ we introduce a new cohomology theory $\,_2H^n(G,A,\kappa;B)$ which we call the secondary cohomology.
We give a construction that associates to a pointed topological space $(X,x_0)$ an  invariant $\,_2\kappa^4\in \,_2H^4(\pi_1(X),\pi_2(X),\kappa^3;\pi_3(X))$. This construction can be seen a ``3-type" generalization of the classical $k$-invariant.

\end{abstract}

\maketitle


%



\section*{Introduction}
The first $k$-invariant $\kappa^3$ of a pointed topological space $(X,x_0)$ was introduced  by Eilenberg and MacLane in \cite{em}. It is a topological invariant that belongs to  $H^3(\pi_1(X), \pi_2(X))$ (where $\pi_n(X)$ is the $n$-th homotopy group). MacLane and Whitehead have proved that equivalence classes of  so called crossed modules are in bijection with the elements of the third cohomology group $H^3(G, A)$. They used this description to show that the $2$-type of a topological space is determined by triples $(\pi_1(X),\pi_2(X),\kappa^3)$ where  $\kappa^3$ the $k$-invariant  mentioned above.

The $3$-type of a space was  characterized by Baues in terms of quadratic modules  (see \cite{hb}).  To our best knowledge there is no description of quadratic modules (or the 3-type) in terms of  cohomology classes. For the classification of $n$-types one has to deal with Postnikov towers and their invariants, which again do not a have purely algebraic description in terms of some cohomology theory.

In this paper we propose a construction that associates to a pointed topological space $(X,x_0)$ an  invariant $\,_2\kappa^4$ that is an element in a certain cohomology group we introduce. The construction is similar with that of the first $k$-invariant  but also has a  Postnikov-invariant flavor. We believe that $\,_2\kappa^4$ is a natural candidate to classify the $3$-type of a space.

Here is how the paper is organized. In the first section we recall general facts about the first $k$-invariant. In the second section, as a warm up, we treat the case of simply connected spaces. More precisely for two commutative  groups $A$ and $B$ we introduce the secondary cohomology  group $_2H^n(A, B)$. Then to a simply connected topological space $X$ we associate a topological invariant $\,_2\kappa^4\in\,_2H^4(\pi_2(X),\pi_3(X))$. This construction is very similar with the construction of the $k$-invariant, one just has to go up one dimension. The key result is the definition of the secondary cohomology groups. Having the right cohomology theory the proof that $\,_2\kappa^4$ is an invariant is almost cut and paste from \cite{em}.

In the third  section we give the result for general  topological spaces.  We start with a group $G$ (possibly noncommutative), two $G$-modules $A$ and $B$, a 3-cocylcle  $\kappa\in H^3(G,A)$ and we define $\,_2H^4(G,A,\kappa;B)$ the secondary cohomology of the triple $(G,A,\kappa)$ with coefficients in $B$. One can then associate to any pointed space $(X,x_0)$ a topological invariant $\,_2\kappa^4\in\,_2H^4(\pi_1(X),\pi_2(X),\kappa^3;\pi_3(X))$. Obviously if $X$ is simply connected we get the invariant $\;_2\kappa^4$ from section two. Also if $\pi_2(X)=0$ we get the second $k$-invariant $\kappa^4\in H^4(\pi_1(X),\pi_3(X))$ introduced in \cite{em}. We conclude the paper with remarks on possible generalizations and research problems in this direction.

\section{Preliminaries}\label{pre}

 We recall from   \cite{em} the construction of the first $k$-invariant. $X$ is an arc connected topological space with base point $x_0$. For each element $\alpha\in \pi_1(X)$ we fix a representative $r(\alpha)$. For each pair of elements $(\alpha,\beta)\in \pi_1(X)\times \pi_1(X)$  we consider a singular 2-simplex $r(\alpha,\beta):\Delta_2\to X$ such that the edges $[0,1]$, $[1,2]$ and $[0,2]$ map according to $r(\alpha)$, $r(\beta)$ and $r(\alpha\beta)$. For $\alpha$, $\beta$, $\gamma \in \pi_1(X)$ we define a map  $R(\alpha,\beta,\gamma):\partial(\Delta_3)\to X$ such that $R_{|[0,1,2]}=r(\alpha,\beta)$, $R_{|[1,2,3]}=r(\beta,\gamma)$, $R_{|[0,2,3]}=r(\alpha\beta,\gamma)$, and $R_{|[0,1,3]}=r(\alpha,\beta\gamma)$. In this way we get an element of $\kappa(\alpha,\beta,\gamma)\in \pi_2(X)$.
\begin{theorem} \cite{em} The cochain $(\alpha,\beta,\gamma)\to \kappa(\alpha,\beta,\gamma)$ is a  cocylcle. A change of the representatives $r(\alpha)$ and $r(\alpha, \beta)$ alters $\kappa$
by a coboundary.  Thus $\kappa$ determines a unique  cohomology class  $\kappa^3\in H^3(\pi_1(X), \pi_2(X))$ which is a topological invariant of $(X,x_0)$.
\end{theorem}

\section{The simply connected case}
\subsection{Secondary cohomology for commutative groups}
In this section $A$ and $B$ are commutative groups. For $A$ we use multiplicative notation while for $B$ we use the additive one.  Define $_2C^n(A,B)=Map(A^{\frac{n(n-1)}{2}},B)$. The elements of  $A^{\frac{n(n-1)}{2}}$ are $\frac{n(n-1)}{2}$-tuples $(a_{i,j})_{(0\leq i<j\leq n-1)}$ with the index in the lexicographic order:
$$(a_{0,1},a_{0,2},...,a_{0,n-1},a_{1,2},a_{1,3},...,a_{1,n-1},...,a_{n-2,n-1})$$
For every $0\leq k\leq n+1$ we define $d_n^k: A^{\frac{(n+1)n}{2}}\to A^{\frac{n(n-1)}{2}}$, $d_n^k((a_{i,j})_{(0\leq i<j\leq n)}=(b_{i,j})_{(0\leq i<j\leq n-1)}$ where
$$b_{i,j}= \left\{
\begin{array}{c l}
a_{i,j}\; {\rm if}\; 0\leq i<j < k-1\\
a_{i,k-1}a_{i,k}a_{k-1,k}^{-1}\; {\rm if}\; 0\leq i<j= k-1\\
a_{i,j+1}\; {\rm if}\; 0\leq i\leq k-1 < j\\
a_{i+1,j+1}\; {\rm if}\; k-1<i <j
\end{array}
\right.
$$
One can check that
\begin{eqnarray*}
d_{n-1}^kd_n^l=d_{n-1}^{l-1}d_n^k \; {\rm if}\; k<l
\end{eqnarray*}
Let $\delta_n:\, _2C^n(A,B)\to\,  _2C^{n+1}(A,B)$ defined by:
\begin{eqnarray}
\delta_n(f)=fd_n^0-fd_n^1+fd_n^2-...+(-1)^{n+1}fd_n^{n+1}
\end{eqnarray}

\begin{example} When $n=2$, $n=3$ or $n=4$, and $f\in \,_2C^n(A,B)$   we have
\begin{eqnarray*}
\delta_2(f)(a_{01},a_{02},a_{12})=f(a_{12})-f(a_{02})
+f(a_{01}a_{02}a_{12}^{-1})-f(a_{01})
\end{eqnarray*}
\begin{eqnarray*}
&&\delta_3(f)(a_{01},a_{02},a_{03},a_{12},a_{13},a_{23})=f(a_{12},a_{13},a_{23})-f(a_{02},a_{03},a_{23})\\
&&+f(a_{01}a_{02}a_{12}^{-1},a_{03},a_{13})
-f(a_{01},a_{02}a_{03}a_{23}^{-1},a_{12}a_{13}a_{23}^{-1})+f(a_{01},a_{02},a_{12})
\end{eqnarray*}
\begin{eqnarray*}
&&\delta_4(f)(a_{01},a_{02},a_{03},a_{04},a_{12},a_{13},a_{14},a_{23},a_{24},a_{34})=\\
&&f(a_{12},a_{13},a_{14},a_{23},a_{24},a_{34})-f(a_{02},a_{03},a_{04},a_{23},a_{24},a_{34})\\
&&+f(a_{01}a_{02}a_{12}^{-1},a_{03},a_{04},a_{13},a_{14},a_{34})\\
&&-f(a_{01},a_{02}a_{03}a_{23}^{-1},a_{04},a_{12}a_{13}a_{23}^{-1},a_{14},a_{24}) \\
&&+f(a_{01},a_{02},a_{03}a_{04}a_{34}^{-1},a_{12},a_{13}a_{14}a_{34}^{-1},a_{23}a_{24}a_{34}^{-1})\\
&&-f(a_{01},a_{02},a_{03},a_{12},a_{13},a_{23})
\end{eqnarray*}
\end{example}

It is obvious that $\delta_{n+1}\delta_n(f)=0$ for all $n\geq 1$ and all $ f\in \, _2C^{n}(A,B)$, which means that we have a complex $(\,_2C^*(A,B),\delta_*)$.
\begin{definition}
We denote the homology of $(\,_2C^*(A,B),\delta_*)$ by $\;_2H^*(A,B)$ and we call it the {\bf secondary cohomology} of the group $A$ with coefficients in $B$.
\end{definition}

\subsection{The secondary $k$-invariant for simply connected spaces}

Let $X$ be a simply connected topological space. For any element $a\in \pi_2(X)$  we fix a map $r(a):\Delta_2\to X$ that represents $a$ (notice that $r(a)_{|\partial(\Delta_2)}=x_0$). For each $a_{01}$, $a_{02}$ and $a_{12}\in \pi_2(X)$ we fix a singular $3$-simplex $r(a_{01},a_{02},a_{12}):\Delta_3\to X$ such that $r(a_{01},a_{02},a_{12})_{|[0,1,2]}=r(a_{01})$, $r(a_{01},a_{02},a_{12})_{|[0,2,3]}=r(a_{02})$,  $r(a_{01},a_{02},a_{12})_{|[1,2,3]}=r(a_{12})$  and $r(a_{01},a_{02},a_{12})_{|[0,1,3]}=r(a_{01}a_{02}a_{12}^{-1})$.
For each $a_{01}$, $a_{02}$, $a_{03}$, $a_{12}$,  $a_{13}$ and  $a_{23}\in \pi_2(X)$  we define:
\begin{eqnarray}
R(a_{01},a_{02},a_{03},a_{12},a_{13},a_{23}):\partial(\Delta_4)\to X
\end{eqnarray}
such that the restriction of $R(a_{01},a_{02},a_{03},a_{12},a_{13},a_{23})$ on each of the five $3$-simplex that make the boundary of $\Delta_4$ is given  by:
\begin{eqnarray*}
&&R(a_{01},a_{02},a_{03},a_{12},a_{13},a_{23})_{|[0,1,2,3]}=r(a_{01},a_{02},a_{12})\\
&&R(a_{01},a_{02},a_{03},a_{12},a_{13},a_{23})_{|[0,1,2,4]}=r(a_{01},
a_{02}a_{03}a_{23}^{-1},a_{12}a_{13}a_{23}^{-1})\\
&&R(a_{01},a_{02},a_{03},a_{12},a_{13},a_{23})_{|[0,1,3,4]}=r(a_{01}a_{02}a_{12}^{-1},a_{03},a_{13})\\
&&R(a_{01},a_{02},a_{03},a_{12},a_{13},a_{23})_{|[0,2,3,4]}=r(a_{02},a_{03},a_{23})\\
&&R(a_{01},a_{02},a_{03},a_{12},a_{13},a_{23})_{|[1,2,3,4]}=r(a_{12},a_{13},a_{23})
\end{eqnarray*}
It is obvious that  $R(a_{01},a_{02},a_{03},a_{12},a_{13},a_{23})$ determines  a unique element of $\pi_3(X)$, and so $R\in \,_2C^4(\pi_2(X),\pi_3(X))$.

Let's  see that $R$ is a cocycle. Take $a_{01}$, $a_{02}$, $a_{03}$, $a_{04}$, $a_{12}$, $a_{13}$, $a_{14}$, $a_{23}$ $a_{24}$ and $a_{34}\in \pi_2(X)$. We notice that there exists a  map $F$ from the $3$-dimensional skeleton  of $\Delta_5$ to $X$ such that:
\begin{eqnarray*}
&&F_{|\partial([0,1,2,3,4])}=R(a_{01},a_{02},a_{03},a_{12},a_{13},a_{23})\\
&&F_{|\partial([0,1,2,3,5])}=R(a_{01},a_{02},a_{03}a_{04}a_{34}^{-1},a_{12},a_{13}a_{14}a_{34}^{-1},a_{23}a_{24}a_{34}^{-1})\\
&&F_{|\partial([0,1,2,4,5])}=R(a_{01},a_{02}a_{03}a_{23}^{-1},a_{04},a_{12}a_{13}a_{23}^{-1},a_{14},a_{24})\\
&&F_{|\partial([0,1,3,4,5])}=R(a_{01}a_{02}a_{12}^{-1},a_{03},a_{04},a_{13},a_{14},a_{34})\\
&&F_{|\partial([0,2,3,4,5])}=R(a_{02},a_{03},a_{04},a_{23},a_{24},a_{34})\\
&&F_{|\partial([1,2,3,4,5])}=R(a_{12},a_{13},a_{14},a_{23},a_{24},a_{34})
\end{eqnarray*}
Moreover each $3$-simplex of $\Delta_5$ appears exactly twice (once for each  orientation) in the following element of $\pi_3(X)$.
\begin{eqnarray*}
&&R(a_{12},a_{13},a_{14},a_{23},a_{24},a_{34})-R(a_{02},a_{03},a_{04},a_{23},a_{24},a_{34})\\
&&+R(a_{01}a_{02}a_{12}^{-1},a_{03},a_{04},a_{13},a_{14},a_{34})\\
&&-R(a_{01},a_{02}a_{03}a_{23}^{-1},a_{04},a_{12}a_{13}a_{23}^{-1},a_{14},a_{24}) \\  &&+R(a_{01},a_{02},a_{03}a_{04}a_{34}^{-1},a_{12},a_{13}a_{14}a_{34}^{-1},a_{23}a_{24}a_{34}^{-1})\\
&&-R(a_{01},a_{02},a_{03},a_{12},a_{13},a_{23})
\end{eqnarray*}
This means that $\delta_4(R)=0$ and so $R\in \,_2C^4(\pi_2(X),\pi_3(X))$ is a $4$-cocylcle.

If we keep fixed $r(a)$ and we change $r(a,b,c)$  with another map $r'(a,b,c)$  we get a map $h:\pi_2(X)\times \pi_2(X)\times\pi_2(X)\to \pi_3(X)$ (by gluing $r$ and $r'$ along the boundary). One can see that:
\begin{eqnarray*}
&&R(a_{01},a_{02},a_{03},a_{12},a_{13},a_{23})-R'(a_{01},a_{02},a_{03},a_{12},a_{13},a_{23})\\
&&=h(a_{12},a_{13},a_{23})-h(a_{02},a_{03},a_{23})+h(a_{01}a_{02}a_{12}^{-1},a_{03},a_{13})\\
&&-h(a_{01},a_{02}a_{03}a_{23}^{-1},a_{12}a_{13}a_{23}^{-1})+
h(a_{01},a_{02},a_{12})
\end{eqnarray*}
And so $R$ and $R'$ are cohomologus equivalent. If we change $r(a)$ with $r'(a)$ then we can chose $r'(a,b,c)$ such that the two maps  $R, R':\pi_2(X)^6\to \pi_3(X)$ are equal. This prove that $R$ defines an unique element $\,_2\kappa^4 \in \,_2H^4(\pi_2(X),\pi_3(X))$ that is a topological invariant of $(X,x_0)$.

\section{The general case}

\subsection{Secondary cohomology of $(G,A,\kappa)$ with coefficients in $B$}

We want to generalize the above results to  topological spaces with $\pi_1(X)$ nontrivial. First we need to construct an analog for $\;_2H^4(A,B)$. We start with a group $G$ (possibly noncommutative), two $G$-modules  $A$ and $B$ and a $3$-cocycle  $\kappa \in Z^3(G,A)$. Define $_2C^n(G,A,\kappa; B)=Map(G^n\times A^{\frac{n(n-1)}{2}},B)$. The elements of $G^n$ are $n$-tuples
$$(g)=(g_i)_{(1\leq i\leq n)}$$
The elements of  $A^{\frac{n(n-1)}{2}}$ are $\frac{n(n-1)}{2}$-tuples $(a)=(a_{i,j})_{(0\leq i<j\leq n-1)}$ with the index in the lexicographic order:
$$(a)=(a_{0,1},a_{0,2},...,a_{0,n-1},a_{1,2},a_{1,3},...,a_{1,n-1},...,a_{n-2,n-1})$$

For every $0\leq k\leq n+1$ we define $d_n^k: G^{n+1}\times A^{\frac{n(n+1)}{2}}\to G^n\times A^{\frac{(n-1)n}{2}}$, $d_n^k((g_i)_{(1\leq i\leq n+1)},(a_{i,j})_{(0\leq i<j\leq n)})=((h_i)_{(1\leq i\leq n)},(b_{i,j})_{(0\leq i<j\leq n-1)})$ where

$$h_{i}= \left\{
\begin{array}{c l}
g_{i}\; {\rm if}\; i < k\\
g_ig_{i+1}\; {\rm if}\; i= k\\
g_{i+1}\; {\rm if}\; k < i\\
\end{array}
\right.
$$
$$b_{i,j}= \left\{
\begin{array}{c l}
a_{i,j}\; {\rm if}\; 0\leq i<j < k-1\\
a_{i,k-1}a_{i,k}\,^{g_{i+1}...g_{k-1}}(a_{k-1,k}^{-1})\kappa(g_{i+1}...g_{k-1},g_k,g_{k+1})\; {\rm if}\; 0\leq i<j= k-1\\
a_{i,j+1}\; {\rm if}\; 0\leq i\leq k-1 < j\\
a_{i+1,j+1}\; {\rm if}\; k-1<i <j
\end{array}
\right.
$$
Let $\delta_n^{\kappa}:\,_2C^n(G,A,\kappa; B)\to\,_2C^{n+1}(G,A,\kappa; B)$ defined by:
\begin{eqnarray*}
\delta_n^{\kappa}(f)((g);(a))=g_1fd_n^0((g);(a))-fd_n^1((g);(a))\\
+fd_n^2((g);(a))-...+(-1)^{n+1}fd_n^{n+1}((g);(a))
\end{eqnarray*}
\begin{example} For $n=2$ or $n=3$ and $f\in \,_2C^n(G,A,\kappa; B)$ we have:
\begin{eqnarray*}
&&\delta_2^{\kappa}(f)(g_1,g_2,g_3;a_{01},a_{02},a_{12})=g_1f(g_2,g_3;a_{12})-f(g_1g_2,g_3;a_{02})\\
&&+f(g_1,g_2g_3;a_{01}a_{02}\,^{g_1}(a_{12}^{-1})\kappa(g_1,g_2,g_3))-f(g_1,g_2;a_{01})
\end{eqnarray*}
\begin{eqnarray*}
&&\delta_3^{\kappa}(f)(g_1,g_2,g_3,g_4;a_{01},a_{02},a_{03},a_{12},a_{13},a_{23})\\
&&=g_1f(g_2,g_3,g_4;a_{12},a_{13},a_{23})-f(g_1g_2,g_3,g_4;a_{02},a_{03},a_{23})\\
&&+f(g_1,g_2g_3,g_4;a_{01}a_{02}\,^{g_1}(a_{12}^{-1})\kappa(g_1,g_2,g_3),a_{03},a_{13})\\
&&-f(g_1,g_2,g_3g_4;a_{01},a_{02}a_{03}\,^{g_1g_2}(a_{23}^{-1})\kappa(g_1g_2,g_3,g_4),a_{12}a_{13}\,^{g_2}(a_{23}^{-1})\kappa(g_2,g_3,g_4))\\
&&+f(g_1,g_2,g_3;a_{01},a_{02},a_{12})
\end{eqnarray*}
\end{example}

One can check that $\delta_{n+1}^{\kappa}\delta_n^{\kappa}(f)=0$ for all $ f\in \, _2C^n(G,A,\kappa; B)$, and so we have a complex $(\,_2C^*(G,A,\kappa; B),\delta_*^{\kappa})$.
\begin{definition} We denote the
homology of the complex $(\,_2C^*(G,A,\kappa; B),\delta_*^{\kappa})$ by $\,_2H^*(G,A,\kappa; B)$ and we call it the {\bf secondary cohomology} of $(G,A,\kappa)$ with coefficients in $B$.
\end{definition}
\begin{remark}
Let's notice  that the above construction depends only on the class of $\kappa \in H^3(G,A)$. Indeed if $\kappa=\kappa'+\delta_2(u)$ then we can define an isomorphism  of complexes $\Phi_u:\,_2C^*(G,A,\kappa; B)\to \,_2C^*(G,A,\kappa'; B)$ defined by:
\begin{eqnarray}
\Phi_u(f)((g);(a))=f((g);(c))\label{nau}
\end{eqnarray}
where $c_{i,j}=a_{i,j}u(g_{i+1}...g_{j-1},g_j)$. One can see that and \begin{eqnarray}
\delta^{\kappa'}\Phi_u=\Phi_u\delta^{\kappa}
\end{eqnarray}
and
\begin{eqnarray}
\Phi_u\Phi_v=\Phi_{u+v}
\end{eqnarray}
And so $\Phi$ is a natural transformation that allows us to identify $\,_2H^*(G,A,\kappa; B)$ with $\,_2H^*(G,A,\kappa'; B)$.
\end{remark}

\begin{example}
If $A$ is trivial then $\,_2H^n(G,1,\kappa; B)$ is the usual cohomology $H^n(G,B)$. If $G$ is trivial then
$\,_2H^*(1,A,\kappa; B)$ is the secondary cohomology group $\,_2H^n(A,B)$ defined in the previous section. Also it is easy to show that $\,_2H^n(1;B)=0$ and $\,_2H^n(A;0)=0$.
\end{example}
\begin{example} Simple computations show that
\begin{equation*}
\,_2H^2(\mathbb{Z}_2,B)=\{(b_1,b_2)\in B\times B\vert \; 2b_1=2b_2\}/\{(b,b)\vert b\in B\}
\end{equation*}
\begin{eqnarray*}
\,_2H^3(\mathbb{Z}_2,B)=B/2B
\end{eqnarray*}
For example one has $\,_2H^2(\mathbb{Z}_2,\mathbb{Z})=0$, $\,_2H^2(\mathbb{Z}_2,\mathbb{Z}_2)=\mathbb{Z}_2$, $\,_2H^3(\mathbb{Z}_2,\mathbb{Z})=\mathbb{Z}_2$ and $\,_2H^3(\mathbb{Z}_2,\mathbb{Z}_2)=\mathbb{Z}_2$.
\end{example}

\subsection{Secondary $k$-invariant} Let $(X,x_0)$ be a  pointed topological space. For each $\alpha\in \pi_1(X)$ we fix a representative $r(\alpha):[0,1]\to X$. For each pair of elements $\alpha$, $\beta\in \pi_1(X)$ we fix a singular $2$-simplex $r(\alpha,\beta):\Delta_2\to X$ such that $[0,1]$, $[1,2]$  and $[0,2]$ map according to $r(\alpha)$, $r(\beta)$ and $r(\alpha\beta)$. Just like in the construction of the $k$-invariant define a map $R(\alpha,\beta,\gamma):\partial(\Delta_3)\to X$ such that  $R_{|[0,1,2]}=r(\alpha,\beta)$, $R_{|[1,2,3]}=r(\beta,\gamma)$, $R_{|[0,2,3]}=r(\alpha\beta,\gamma)$, and $R_{|[0,1,3]}=r(\alpha,\beta\gamma)$. This gives us the classical $k$-invariant.

For each triple $(\alpha,\beta;a)\in \pi_1(X)\times \pi_1(X)\times \pi_2(X)$ we consider a singular $2$-simplex $r(\alpha,\beta;a):\Delta_2\to X$ such that $[0,1]$, $[1,2]$  and $[0,2]$ map according to $r(\alpha)$, $r(\beta)$ and $r(\alpha\beta)$ and when we glue $r(\alpha,\beta;a)$ with $r(\alpha,\beta)$ along the boundary we get $a\in \pi_2(X)$. For each $(g)=(g_1,g_2,g_3)\in \pi_1(X)^3$  and $(a)=(a_{01},a_{02},a_{12})\in \pi_2(X)^3$  we fix a singular $3$-simplex
$r((g);(a))=r(g_1,g_2,g_3;a_{01},a_{02},a_{12}):\Delta_3\to X$  such that:
\begin{eqnarray*}
r((g);(a))_{|[0,1,2]}&=&r(g_1,g_2;a_{01})\\ r((g);(a))_{|[0,2,3]}&=&r(g_1g_2,g_3;a_{02})\\  r((g);(a))_{|[1,2,3]}&=&r(g_2,g_3;a_{12})\\   r((g);(a))_{|[0,1,3]}&=&r(g_1,g_2g_3;a_{01}a_{02}\,^{g_1}(a_{12}^{-1})\kappa(g_1,g_2,g_3))
\end{eqnarray*}

For each $(g)=(g_1,g_2,g_3,g_4)\in \pi_1(X)^4$ and $(a)=(a_{01}, a_{02}, a_{03}, a_{12},  a_{13}, a_{23})\in \pi_2(X)^6$  we define:
\begin{eqnarray}
R((g);(a))=R(g_1,g_2,g_3,g_4;a_{01},a_{02},a_{03},a_{12},a_{13},a_{23}):\partial(\Delta_4)\to X
\end{eqnarray}
such that the restriction of $R((g);(a))$ on each of the five $3$-simplex that make the boundary of $\Delta_4$ is given  by:
\begin{eqnarray*}
&&R((g);(a))_{|[0,1,2,3]}=r(g_1,g_2,g_3;a_{01},a_{02},a_{12})\\
&&R((g);(a))_{|[0,1,2,4]}=r(g_1,g_2,g_3g_4;a_{01},
a_{02}a_{03}\,^{g_1g_2}(a_{23}^{-1})\kappa(g_1g_2,g_3,g_4),\\
&&a_{12}a_{13}\,^{g_2}(a_{23}^{-1})\kappa(g_2,g_3,g_4))\\
&&R((g);(a))_{|[0,1,3,4]}=r(g_1,g_2g_3,g_4;a_{01}a_{02}\,^{g_1}(a_{12}^{-1})\kappa(g_1,g_2,g_3),a_{03},a_{13})\\
&&R((g);(a))_{|[0,2,3,4]}=r(g_1g_2,g_3,g_4;a_{02},a_{03},a_{23})\\
&&R((g);(a))_{|[1,2,3,4]}=r(g_2,g_3,g_4;a_{12},a_{13},a_{23})
\end{eqnarray*}
It is obvious that  $R(g_1,g_2,g_3,g_4;a_{01},a_{02},a_{03},a_{12},a_{13},a_{23})$ determines  a unique element of $\pi_3(X)$, and so $R\in \,_2C^4(\pi_1(X),\pi_2(X),\kappa;\pi_3(X))$.

Just like in the case $X$ simply connected one can show that
$R$ is a 4-cocycle, i.e. $R\in \,_2Z^4(\pi_1(X),\pi_2(X),\kappa;\pi_3(X))$.

We want to show that the class of $R\in \,_2H^4(\pi_1(X),\pi_2(X),\kappa;\pi_3(X))$ does not depend on the choices we made. First we keep fixed $r(\alpha)$ and $r(\alpha, \beta)$ and $r(\alpha, \beta; a)$ and we change $r(\alpha,\beta,\gamma; a,b,c)$  with another map $r'(\alpha,\beta,\gamma; a,b,c)$  we get a map $h:\pi_1(X)^3\times \pi_2(X)^3\to \pi_3(X)$ (by gluing $r$ and $r'$ along the boundary). One can see that:
\begin{eqnarray*}
&&R-R'=\delta_3^{\kappa}(h)
\end{eqnarray*}
And so $R$ and $R'$ are cohomologus equivalent. If we change $r(\alpha,\beta; a)$ with $r'(\alpha,\beta; a)$ then we can chose $r'(\alpha,\beta,\gamma;a,b,c)$ such that the two maps  $R, R':\pi_2(X)^6\to \pi_3(X)$ are equal. If we change  $r(\alpha, \beta)$ with $r'(\alpha,\beta))$ we replace  $\kappa$ with some $\kappa'=\kappa-\delta_2(u)$ which gives an isomorphism like in (\ref{nau}). If we fix  $r(\alpha,\beta;a)$ and  $r(\alpha,\beta,\gamma;a,b,c)$ then the 4-cocycle $R'\in \,_2Z^4(\pi_1(X),\pi_2(X),\kappa';\pi_3(X))$ becomes $\Phi_u(R)$. Finally if we change $r(\alpha)$  with $r'(\alpha)$ we can chose all the other $r'$ such that $R$ does not change.   We have the following result.
\begin{theorem} If $(X,x_0)$ is a topological space then the above construction defines a topological invariant  $\,_2\kappa^4 =R\in \,_2H^4(\pi_1(X),\pi_2(X),\kappa;\pi_3(X))$.
\end{theorem}

\begin{remark} If $\pi_2(X)=0$ then $\,_2\kappa^4$ is the element $\kappa^4$ described in the first remark from section 5 in \cite{em} (see also \cite{em2}).
\end{remark}

\subsection{Conclusions and Remarks}

A natural question is whether the invariant  $\,_2\kappa^4$  classify the $3$-type of a space. A possible approach to this problem would be to use the results proved in \cite{hb} and then show that
equivalences classes of quadratic modules are in bijection with elements of the secondary cohomology group $\,_2H^4(G,A,\kappa;B)$.

If the the above question has a positive answer one could try to define a ternary cohomology group $\,_3H^n(G,A,\kappa,B,\,_2\kappa; C)$.  Then  find a cohomology class   $\;_3\kappa^5\in \,_3H^5(\pi_1(X),\pi_2(X),\kappa^3,\pi_3(X),\,_2\kappa^4; \pi_4(X))$ that classify the $4$-type of a space, and so on. We can notice that we have  a short exact sequence of complexes:
\begin{eqnarray*}
0\to \,C^*(G,B)\to \,_2C^*(G,A,\kappa;B)\to \,_2C^*(A,B)\to 0
\end{eqnarray*}
This suggest that at the next level we should have:
\begin{eqnarray*}
0\to \,_2C^*(G,A,\kappa;C)\to \,_3C^*(G,A,\kappa,B,\,_2\kappa; C)
\to \,_3C^*(B, C)\to 0
\end{eqnarray*}
In general we expect that the cohomology theory at step $n$ is a twist between  the cohmology from step $n-1$ with an appropriate cohomology theory that depends only on two groups.  This is  similar with the results from \cite{may} and is also the reason way we said in introduction that our construction has a Postnikov-invariant flavor.

Even if the above questions  do not have a positive answer one could try to to give an algebraic description of $H^3(X,K^*)$ in terms $\pi_1(X)$, $\pi_2(X)$, $\pi_3(X)$, $\kappa$ and $\;_2\kappa^4$, just like $H^2(X,K^*)$ was determined in \cite{em} in terms of $\pi_1(X)$, $\pi_2(X)$ and $\kappa^3$.

Finally, notice that when we prove $\delta_4\delta_3=0$ we use an equality of the type
\begin{eqnarray}
f(f(a_{01},a_{0,2},a_{12}),a_{03},a_{13})=f(a_{01},f(a_{02},a_{03},a_{23}),f(a_{12},a_{13},a_{23})) \label{tas}
\end{eqnarray}
where $f:A\times A\times A\to A$, $f(a,b,c)=abc^{-1}$. This is almost the same as the ternary associativity condition discussed in \cite{sm}. With the notations from that paper one can take  $f(a,b,c)=m(a,Q(c),b)$ and check that $f$ satisfy condition (\ref{tas})




\bibliographystyle{amsalpha}

\begin{thebibliography}{A}

\bibitem
[B]
{hb} H. J. Baues,
\textit{Combinatorial Homotopy and $4$-Dimensional Complexes.}  Walter de Gruyter, Berlin, (1991).

\bibitem
[EM]
{em} S. Eilenberg and S. MacLane
\textit{Determinationation of the Second Homology and Cohomology Groups of a Space by Means of Homotopy Invariants.} Proc. National Academy of Science, vol {\bf 32} (1946) 277-280.

\bibitem
[EM2]
{em2} S. Eilenberg and S. MacLane
\textit{ Relations between homology and homotopy groups of spaces. II.} Ann. of Math. (3) {\bf 51} (1950) 514-533.

\bibitem
[MW]
{mw} S. MacLane and J. H. C. Whitehead \textit{On the 3-type of a Complex}. Proc. National Academy of Science, vol {\bf 36} (1950) 155-178.


\bibitem
[M]
{may} J. P. May,
\textit{Simplicial Objects in Algebraic Topology.} Chicago Lectures in Mathematics (1967).

\bibitem
[SM]
{sm} M. D. Staic,
\textit{From 3-algebras to $\Delta$-groups and Symmetric Cohomology.} Journal of Algebra (4) vol. {\bf 322} (2009) 1360-1378.


\end{thebibliography}

\end{document}